\documentclass[12pt]{article}
\usepackage[cp1251]{inputenc}
\usepackage[russian]{babel}
\usepackage{amssymb}
\usepackage{amsmath}

\begin{document}

\begin{center}
\textbf{\large Absolutely convergent Fourier series. An
improvement of the Beurling--Helson theorem}
\end{center}

\begin{center}
Vladimir Lebedev
\end{center}

\begin{quotation}
{\small \textsc{Abstract.} We consider the space $A(\mathbb T)$ of
all continuous functions $f$ on the circle $\mathbb T$ such that
the sequence of Fourier coefficients
$\widehat{f}=\{\widehat{f}(k), ~k \in \mathbb Z\}$ belongs to
$l^1(\mathbb Z)$. The norm on $A(\mathbb T)$ is defined by
$\|f\|_{A(\mathbb T)}=\|\widehat{f}\|_{l^1(\mathbb Z)}$. According
to the known Beurling--Helson theorem, if $\varphi : \mathbb
T\rightarrow\mathbb T$ is a continuous mapping such that
$\|e^{in\varphi}\|_{A(\mathbb T)}=O(1), ~n\in\mathbb Z,$ then
$\varphi$ is linear. It was conjectured by Kahane that the same
conclusion about $\varphi$ is true under the assumption that
$\|e^{in\varphi}\|_{A(\mathbb T)}=o(\log |n|)$. We show that if
$\|e^{in\varphi}\|_{A(\mathbb T)}=o((\log\log |n|/\log\log\log
|n|)^{1/12})$ then $\varphi$ is linear.

  References: 15 items.


  }
\end{quotation}

\quad

   We consider the class $A(\mathbb T)$ of all continuous functions
$f$ on the circle $\mathbb T=\mathbb R/2\pi\mathbb Z$ whose
Fourier series
$$
f(t)\sim \sum_{k\in\mathbb Z}\widehat{f}(k)e^{ikt}
$$
converge absolutely. Here $\mathbb R$ is the real line, $\mathbb
Z$ is the group of integers and $\widehat{f}=\{\widehat{f}(k),
~k\in\mathbb Z\}$ is the sequence of Fourier coefficients of a
function $f$,
$$
\widehat{f}(k)=\int_{\mathbb T} f(t)e^{-ikt}\frac{dt}{2\pi},
\qquad k\in\mathbb Z.
$$
The class $A(\mathbb T)$ is a Banach space with respect to the
natural norm
$$
\|f\|_{A(\mathbb T)}=\|\widehat{f}\|_{l^1(\mathbb Z)}=
\sum_{k\in\mathbb Z}|\widehat{f}(k)|.
$$
It is well known that $A(\mathbb T)$ is a Banach algebra (with the
usual multiplication of functions).

  We identify continuous mappings $\varphi$ of the circle
$\mathbb T$ into itself and continuous functions $\varphi :
\mathbb R\rightarrow\mathbb R$ satisfying
$$
\varphi(t+2\pi)=\varphi(t)~(\mathrm{mod}~2\pi).
\eqno(1)
$$

  Recall the known Beurling--Helson theorem
[1] (see also [2], [3]). Let $\varphi : \mathbb
T\rightarrow\mathbb T$ be a continuous mapping such that
$$
\|e^{in\varphi}\|_{A(\mathbb T)}=O(1),
\qquad n\in\mathbb Z,
$$
then $\varphi$ is linear (with integer slope) i.e. $\varphi(t)=\nu
t+\varphi(0)$, where $\nu\in\mathbb Z$. This theorem immediately
implies the solution of the Levy problem on the description of
endomorphisms of the algebra $A(\mathbb T)$: all these
endomorphisms are trivial, i.e., have the form $f(t)\rightarrow
f(\nu t+t_0)$. In other words only trivial changes of variable are
admissible in $A(\mathbb T)$. Note also another version of the
Beurling--Helson theorem: if $U$ is a bounded translation
invariant operator from $l^1(\mathbb Z)$ to itself such that
$\|U^n\|_{l^1\rightarrow l^1}=O(1), ~n \in \mathbb Z$, then $U=\xi
S$, where $\xi$ is a complex number, $|\xi|=1$, and $S$ is a
translation.

  The character of the growth of the norms
$\|e^{in\varphi}\|_{A(\mathbb T)}$ is in general not clear. We
shall briefly indicate certain known results. It is easy to show
(see [2, Ch. VI, \S~3]) that for every $C^1$ -smooth mapping
$\varphi$ we have $\|e^{in\varphi}\|_{A(\mathbb
T)}=O(\sqrt{|n|})$. At the same time if a mapping $\varphi$ is
$C^2$ -smooth and nonlinear then
$$
\|e^{in\varphi}\|_{A(\mathbb T)}\geq c\sqrt{|n|},
\qquad n\in\mathbb Z,
$$
where $c=c(\varphi)$ is independent of $n$. This result is
contained implicitly in Z. L. Leibenson's work [4] and was
obtained in explicit form by J. -P. Kahane [5] with the use of
Leibenson's approach (for a simple proof see [2, Ch. VI, \S~3]).
Thus, for every nonlinear $C^2$ -smooth mapping $\varphi$ we have
$$
\|e^{in\varphi}\|_{A(\mathbb T)}\simeq\sqrt{|n|},
\qquad |n|\rightarrow\infty.
$$

  The growth of the norms $\|e^{in\varphi}\|_{A(\mathbb T)}$ in
the  $C^1$ -smooth case was studied by the author of the present
paper in [6], [7] (see also [8] where we considered functions of
several variables). Note that the $C^1$ -smooth case is
essentially different from the $C^2$ -smooth case.

  In general, the norms $\|e^{in\varphi}\|_{A(\mathbb T)}$ can
grow rather slowly. Kahane [2, Ch. VI, \S~2] showed that if
$\varphi: \mathbb T\rightarrow\mathbb T$ is a piecewise linear but
not linear continuous mapping, then $\|e^{in\varphi}\|_{A(\mathbb
T)}\simeq \log |n|$.

   It is unknown if for nonlinear continuous mappings
$\varphi$ the norms $\|e^{in\varphi}\|_{A(\mathbb T)}$ can grow
slower than $\log |n|$. Kahane conjectured that the
Beurling--Helson theorem can be strengthened, in particular, he
conjectured that the condition
$$
\|e^{in\varphi}\|_{A(\mathbb T)}=o(\log
|n|), \quad |n|\rightarrow\infty,
$$
implies that $\varphi$ is linear. As far as the author knows this
conjecture was first proposed by Kahane in his talk presented at
the ICM 1962 in Stockholm, [9]. Later the conjecture was mentioned
by Kahane in [2] and [3].\footnote{Note that a priori the mere
existence of a sequence $\omega_n\rightarrow\infty$ such that the
condition $\|e^{in\varphi}\|_{A(\mathbb T)}=o(\omega_{|n|})$
implies linearity of a mapping $\varphi$ is not obvious.}

  Here we shall obtain the following theorem.

\quad

\textbf{Theorem.} \emph{Let $\varphi: \mathbb T\rightarrow\mathbb
T$ be a continuous mapping. Suppose that
$$
\|e^{in\varphi}\|_{A(\mathbb
T)}=o\bigg(\bigg(\frac{\log\log |n|}{\log\log\log |n|}\bigg)^{1/12}\bigg),
\qquad n\in\mathbb Z, \quad |n|\rightarrow \infty.
$$
Then $\varphi$ is linear, i.e., $\varphi(t)=\nu t+\varphi(0),
~\nu\in\mathbb Z$.}

\quad

  Ideologically the proof of our theorem is to some extent close
to the proof of the Beurling--Helson theorem given by Kahane in
[3] (the proof in [3] is based on completely different approach
rather than the original proof by Beurling and Helson [1], [2]).
We modify Kahane's arguments so that instead of dealing with the
group $\mathbb T$ and the mapping $\varphi$ itself we apply the
arguments to the cyclic group $\mathbb T_N$ for large $N's$ and a
certain mapping $\varphi_N$, which on $\mathbb T_N$ provides a
good approximation for $\varphi$ and whose values are rational
numbers with a ``small common denominator''. The construction of
this mapping is based on the Dirichlet theorem on simultaneous
diophantine approximation. The key role in the proof is played by
the result recently obtained by B. Green and S. V. Konyagin [10,
Theorem 1.3].

    Before proving the theorem let us introduce notation,
recall certain (quite standard, see, e.g., [11]) facts of harmonic
analysis on $\mathbb T_N$, and discuss the Green--Konyagin
theorem.

   Let $G$ be a finite abelian group.  We
use the normalised counting measure $\mu_G$ on the group $G$ which
attaches weight $1/(\mathrm{card} G)$ to each point $x\in G$ (by
$\mathrm{card} X$ we denote the number of the elements of a finite
set $X$). Thus, if $E\subseteq G$, then
$$
\mu_G (E)=\frac{\mathrm{card} E}{\mathrm{card} G}
$$
and if $f$ is a (complex) function on $G$, then
$$
\int_G f(x)d\mu_G(x)=\frac{1}{\mathrm{card} G}\sum_{x\in G} f(x).
$$
We shall often write $\int_G f(x)dx$ or $\int_G fd\mu_G$ instead
of $\int_G f(x)d\mu_G(x)$.

  We put
$$
\|f\|_{L^2(G)}=\bigg(\int_G|f(x)|^2 dx\bigg)^{1/2}, \qquad
\|f\|_{L^\infty(G)}=\max_{x\in G}|f(x)|.
$$
In general, for a set $E\subseteq G$ and a function $f$ on $E$ we
put
$$
\|f\|_{L^\infty(E)}=\max_{x\in E}|f(x)|.
$$

   Any finite abelian group $G$ endowed with the measure $\mu_G$ is
a probability space. It will be convenient to consider the general
case. Let $(\mathbb X,\mu)$ be a probability space. For a set
$E\subseteq \mathbb X$ we denote by $1_E$ its characteristic
function: $1_E(x)=1$ if $x\in E$ and $1_E(x)=0$ if $x\notin E$.
For an arbitrary $\mu$ -measurable set $E\subseteq\mathbb X$ let
$$
\delta_{\mathbb X, \mu}(E)=\min(\mu(E), 1-\mu(E)).
$$
Clearly we always have $0\leq\delta_{\mathbb X, \mu}(E)\leq1/2$.
If $\delta_{\mathbb X, \mu}(E)$ is very small then the set $E$ is
either ``very small'' or ``very large'' (i.e., its compliment is
small).

   We denote by $\mathbb T_N$ the cyclic group
$$
\mathbb T_N=\bigg\{\frac{2\pi j}{N}, ~j=0, 1, \ldots, N-1\bigg\}
$$
regarded as a subgroup of $\mathbb T$ (with addition
$\mathrm{mod}~2\pi$). The characters of $\mathbb T_N$ are the
functions $e_k, ~k=0, 1, \ldots, N-1,$ defined as follows:
$$
e_k\bigg(\frac{2\pi j}{N}\bigg)=e^{ik\frac{2\pi j}{N}},
\qquad j=0, 1, \ldots, N-1.
$$
We realize the group dual to $\mathbb T_N$ as the group
$$
\mathbb Z_N=\{0, 1, \ldots, N-1\}
$$
considered with addition $\mathrm{mod}~N$. If $f$ is a function on
$\mathbb T_N$, then its Fourier transform
$\widehat{f}=\{\widehat{f}(k), ~k\in\mathbb Z_N\}$ on the group
$\mathbb T_N$ is defined by
$$
\widehat{f}(k)=\int_{\mathbb T_N} f \overline{e_k} d\mu_{\mathbb T_N},
\qquad k\in\mathbb Z_N
$$
(the bar stands for the complex conjugation). We use the same
notation $\widehat{f}$ for the Fourier transform on $\mathbb T$
and on $\mathbb T_N$ but this does not lead to confusion; in what
follows, up to the end of the proof of the theorem, $\widehat{f}$
stands for the Fourier transform on $\mathbb T_N$.

  Let $f$ be an arbitrary (complex) function on
$\mathbb T_N$. We have
$$
f(t)=\sum_{k\in\mathbb Z_N}\widehat{f}(k)e_k(t),
\qquad t\in \mathbb T_N.
$$
Put
$$
\|f\|_{A(\mathbb T_N)}=\|\widehat{f}\|_{l^1(\mathbb Z_N)}=
\sum_{k\in\mathbb Z_N}|\widehat{f}(k)|.
$$
It is easy to verify that for any two functions $f_1, ~f_2$ on
$\mathbb T_N$ we have
$$
\|f_1f_2\|_{A(\mathbb
T_N)}\leq\|f_1\|_{A(\mathbb T_N)}\|f_2\|_{A(\mathbb T_N)}.
$$

   Generally, for $p\geq 1$ we put
$$
\|f\|_{A_p(\mathbb T_N)}=\|\widehat{f}\|_{l^p(\mathbb Z_N)}=
\bigg(\sum_{k\in\mathbb Z_N}|\widehat{f}(k)|^p\bigg)^{1/p}
$$
(we write $A(\mathbb T_N)$ instead of $A_1(\mathbb T_N)$).

  The Parseval identity has the form
$$
\|f\|_{L^2(\mathbb T_N)}=\|f\|_{A_2(\mathbb T_N)}.
$$

  For any two functions $f_1, f_2$ on $\mathbb T_N$ we
define their convolution $f_1\ast f_2$ by
$$
f_1\ast f_2 (t)=\int_{\mathbb T_N} f_1(x)f_2(t-x)dx,
\qquad t\in \mathbb T_N.
$$
We have $\widehat{f_1\ast f_2}=\widehat{f_1}\widehat{f_2}$.

  Now we recall the Green--Konyagin theorem. Note first the
trivial fact: if $E$ is a set on the circle $\mathbb T$ such that
the Fourier transform of its characteristic function $1_E$ belongs
to $l^1(\mathbb Z)$, then either the set $E$ has (Lebesgue)
measure zero or it is almost the whole circle $\mathbb T$. This
effect can be viewed in more detail on the cyclic group $\mathbb
T_N$. Namely, according to the Green--Konyagin theorem, if $N$ is
a sufficiently large prime number, then for any real function $f$
on $\mathbb T_N$ with $\int_{\mathbb T_N}f d\mu_{\mathbb T_N}=0$
we have
$$
\min_{t\in\mathbb T_N}|f(t)|\leq
c\bigg(\frac{\log\log N}{\log N}\bigg)^{1/3}\|f\|_{A(\mathbb T_N)},
$$
where the constant $c>0$ is independent of $N$ and $f$. Taking an
arbitrary set $E\subseteq\mathbb T_N$ and applying this estimate
to $f(t)=1_E(t)-\mu_{\mathbb T_N}(E)$, we see that
$$
\delta_{\mathbb T_N, ~\mu_{\mathbb T_N}}(E)\leq
c \bigg(\frac{\log\log N}{\log N}\bigg)^{1/3}\|1_E\|_{A(\mathbb T_N)},
\eqno(2)
$$
where $c>0$ is independent of $N$ and $E$. (We took into account
that $\mu_{\mathbb T_N}(E)=\widehat{1_E}(0)$ whence
$\|f\|_{A(\mathbb T_N)}\leq\|1_E\|_{A(\mathbb T_N)}$.)

  We shall also need the Dirichlet theorem on
simultaneous diophantine approximation (see, e.g., [12, Ch. II,
\S~1]). If $\alpha_1, \alpha_2, \ldots, \alpha_N$ are real numbers
and $D\geq 1$ is an integer, then there exist integers $Q, P_1,
P_2, \ldots, P_N$ such that $1\leq Q\leq D^N$ and
$$
|\alpha_jQ-P_j|\leq\frac{1}{D}, \qquad j=1, 2, \ldots, N.
$$

\quad

  \emph{Proof of the Theorem.} Since the mapping $\varphi$ is
continuous, we have
$$
\varphi(t+2\pi)=\varphi(t)+2\pi k,
$$
where $k\in\mathbb Z$ is independent of $t$. Replacing the
function $\varphi$ by $\varphi(t)-kt$, we can assume that
$$
\varphi(t+2\pi)=\varphi(t),
$$
i.e., that $\varphi$ is a real continuous function on the circle
$\mathbb T$.

   Note also that just by the assumption of the theorem
we have $e^{in\varphi}\in A(\mathbb T)$ for all sufficiently large
$n\in\mathbb Z$. Taking sufficiently large $n$ we have
$e^{in\varphi}\in A(\mathbb T)$ and $e^{i(n+1)\varphi}\in
A(\mathbb T)$. Then $e^{-in\varphi}=\overline{e^{in\varphi}}\in
A(\mathbb T)$, and we see that
$e^{i\varphi}=e^{i(n+1)\varphi}e^{-in\varphi}\in A(\mathbb T)$.
Therefore, $e^{in\varphi}\in A(\mathbb T)$ for all $n\in\mathbb
Z$.

  We put
$$
\Theta(n)=\max_{\lambda=0, 1, 2, \ldots, n}
\|e^{i\lambda\varphi}\|_{A(\mathbb T)}.
$$
The sequence $\Theta(n), ~n=0, 1, 2, \ldots,$ is non-decreasing,
$\Theta(n)\geq\Theta(0)=1$ for all $n$, and, by the assumption of
the theorem,
$$
\Theta(n)=
o\bigg(\bigg(\frac{\log\log n}{\log\log\log n}\bigg)^{1/12}\bigg).
\eqno(3)
$$

    We have
$$
\|e^{in\varphi}\|_{A(\mathbb T)}\leq \Theta(n),
\quad n=0, 1, 2, \ldots,
\eqno(4)
$$

  Everywhere below in the proof up to Lemma 4 we assume that
$N$ is an arbitrary positive integer.

  Fixing $N$ and applying the Dirichlet theorem to the numbers
$$
\alpha_j=\frac{1}{2\pi}\varphi\bigg(\frac{2\pi j}{N}\bigg),
\quad j=0, 1, \ldots, N-1,
$$
and $D=N$, we find integers $Q_N$ and $P_{N, j}, ~j=0, 1, \ldots,
N-1,$ such that
$$
1\leq Q_N\leq N^N
\eqno(5)
$$
and
$$
\bigg|\frac{1}{2\pi}\varphi\bigg(\frac{2\pi j}{N}\bigg)Q_N-P_{N, j}\bigg|
\leq\frac{1}{N}, \qquad j=0, 1, \ldots, N-1.
$$

   Let us define a function $\varphi_N$ on $\mathbb T_N$ by
$$
\varphi_N\bigg(\frac{2\pi j}{N}\bigg)=\frac{2\pi P_{N, j}}{Q_N},
\qquad j=0, 1, \ldots, N-1.
\eqno(6)
$$
We see that
$$
\bigg|\varphi\bigg(\frac{2\pi j}{N}\bigg)-
\varphi_N\bigg(\frac{2\pi j}{N}\bigg)\bigg|\leq\frac{2\pi}{NQ_N},
\qquad j=0, 1, \ldots, N-1.
$$
That is
$$
\|\varphi-\varphi_N\|_{L^\infty(\mathbb T_N)}\leq\frac{2\pi}{NQ_N}.
\eqno(7)
$$

\quad

\textbf{Lemma 1.} \emph{For all $n=0,1,2, \ldots, Q_N-1$ we have
$\|e^{in\varphi_N}\|_{A(\mathbb T_N)}\leq 8 \Theta(N^N)$.}

\quad

\emph{Proof.} Let $f$ be an arbitrary function in $A(\mathbb T)$.
We have
$$
f(t)=\sum_{k\in\mathbb Z}c_k e^{ikt},
\qquad \sum_{k\in\mathbb Z} |c_k|=\|f\|_{A(\mathbb T)}<\infty.
$$
Considering $f$ as a function on $\mathbb T_N$ and calculating its
Fourier transform $\widehat{f}$ on $\mathbb T_N$, we obtain
\begin{multline*}
\widehat{f}(k)= \int_{\mathbb T_N}\bigg(\sum_{\nu\in\mathbb Z}
c_\nu e^{i\nu t}\bigg)e^{-ikt} d\mu_{\mathbb T_N}(t)\\=
\sum_{\nu\in\mathbb Z}c_\nu\int_{\mathbb T_N}e^{i\nu t}e^{-ikt}
d\mu_{\mathbb T_N}(t)=\sum_{\nu=k(\mathrm{mod} N)}c_\nu, \qquad
k\in\mathbb Z_N,
\end{multline*}
therefore,
$$
\|f\|_{A(\mathbb T_N)}\leq\|f\|_{A(\mathbb T)}.
$$
Using this inequality, we see that (see (4))
$$
\|e^{in\varphi}\|_{A(\mathbb T_N)}\leq \Theta(n),
\qquad n=0, 1, 2, \ldots.
\eqno(8)
$$

  Note that estimate (7) yields (for $n\geq 0$)
$$
\|e^{in\varphi}-e^{in\varphi_N}\|_{L^\infty(\mathbb T_N)}
\leq n\|\varphi-\varphi_N\|_{L^\infty(\mathbb T_N)}
\leq n\frac{2\pi}{NQ_N},
$$
whence for $n=0, 1, \ldots, Q_N-1$ we obtain
$$
\|e^{in\varphi}-e^{in\varphi_N}\|_{L^\infty(\mathbb T_N)}
\leq\frac{2\pi}{N}.
\eqno(9)
$$

  Clearly, for each function $f$ on $\mathbb T_N$ we have
$$
|\widehat{f}(k)|\leq
\|f\|_{L^\infty(\mathbb T_N)}, \qquad k\in\mathbb Z_N,
$$
($\widehat{f}$ is the Fourier transform on $\mathbb T_N$) whence,
taking into account that $\mathrm{card}~\mathbb Z_N=N$, we see
that
$$
\|f\|_{A(\mathbb T_N)}=\sum_{k\in\mathbb Z_N}|\widehat{f}(k)|
\leq N \|f\|_{L^\infty(\mathbb T_N)}.
$$
Therefore (see (9)),
$$
\|e^{in\varphi}-e^{in\varphi_N}\|_{A(\mathbb T_N)}
\leq 2\pi, \qquad n=0, 1, \ldots, Q_N-1.
\eqno(10)
$$

   Thus, for $n=0, 1, \ldots, Q_N-1$ we obtain
(see (8), (10) and (5))
\begin{multline*}
\|e^{in\varphi_N}\|_{A(\mathbb T_N)}\leq
\|e^{in\varphi}\|_{A(\mathbb T_N)}+
\|e^{in\varphi}-e^{in\varphi_N}\|_{A(\mathbb T_N)} \\ \leq
\Theta(n)+2\pi\leq 8\Theta(n)\leq 8\Theta(Q_N)\leq 8\Theta(N^N).
\end{multline*}
The lemma is proved.

\quad

  We define a function $\Phi$ on
$\mathbb T^3=\mathbb T\times\mathbb T\times\mathbb T$ by
$$
\Phi(x, y, z)=\varphi(x)+\varphi(z-x)-\varphi(y)-\varphi(z-y),
\qquad x, y, z\in \mathbb T.
$$

   We also define a function $\Phi_N$ on the group
$\mathbb T_N^3=\mathbb T_N\times\mathbb T_N\times\mathbb T_N$, by
$$
\Phi_N(x, y, z)=\varphi_N(x)+\varphi_N(z-x)-\varphi_N(y)-\varphi_N(z-y),
\qquad x, y, z \in\mathbb T_N.
$$

  Consider the set
$$
E_N=\{(x, y, z)\in\mathbb T_N^3 : e^{i\Phi_N(x, y, z)}=1\}.
$$

   Note that each value of the function $\varphi_N$ and thus
each value of the function $\Phi_N$ as well is of the form $2\pi
P/Q_N$ for a certain $P\in\mathbb Z$ (see (6)), so the following
identity holds:
$$
\frac{1}{Q_N}\sum_{n=0}^{Q_N-1}e^{in\Phi_N}=1_{E_N}.
\eqno(11)
$$

  In the following lemma we establish the lower bound
for the measure of the set $E_N$.

\quad

\textbf{Lemma 2.} \emph{We have}
$$
\frac{1}{64(\Theta(N^N))^2}\leq \mu_{\mathbb T_N^3} (E_N).
$$

\quad

\emph{Proof.} First let us verify that if $f$ is an arbitrary real
function on $\mathbb T_N$, then
$$
\frac{1}{\|e^{inf}\|_{A(\mathbb T_N)}^2}\leq
\int_{\mathbb T_N^3} e^{inF} d\mu_{\mathbb T_N^3},
\qquad n\in\mathbb Z,
\eqno(12)
$$
where
$$
F(x, y, z)=f(x)+f(z-x)-f(y)-f(z-y),
\qquad x, y, z\in \mathbb T_N.
$$
To see this we repeat word for word Kahane's arguments [3] just
replacing $\mathbb T$ with $\mathbb T_N$. Namely, interpolating
$l^2$ between $l^1$ and $l^4$, we obtain $\|\cdot\|_{l^2}\leq
\|\cdot\|_{l^1}^{1/3}\|\cdot\|_{l^4}^{2/3}$. Therefore,
$$
1=\|e^{inf}\|_{L^2(\mathbb T_N)}=\|e^{inf}\|_{A_2(\mathbb T_N)}
\leq\|e^{inf}\|_{A(\mathbb T_N)}^{1/3}
\|e^{inf}\|_{A_4(\mathbb T_N)}^{2/3}.
$$
Hence
\begin{multline*}
\frac{1}{\|e^{inf}\|_{A(\mathbb T_N)}^2}\leq
\|e^{inf}\|_{A_4(\mathbb T_N)}^4= \sum_{k\in \mathbb
Z_N}|\widehat{e^{inf}\ast e^{inf}}(k)|^2\\= \|e^{inf}\ast
e^{inf}\|_{L^2(\mathbb T_N)}^2= \int_{\mathbb T_N}
\bigg|\int_{\mathbb T_N} e^{inf(x)}e^{inf(z-x)} dx\bigg|^2 dz
\\ = \int_{\mathbb T_N}\bigg(\int_{\mathbb T_N}
e^{inf(x)}e^{inf(z-x)} dx\bigg)\overline{\bigg(\int_{\mathbb T_N}
e^{inf(y)}e^{inf(z-y)} dy\bigg)} dz \\ = \iiint_{\mathbb T_N^3}
e^{inF(x, y, z)} dx dy dz.
\end{multline*}
Thus we have (12).

    Now we use inequality (12) for
$f=\varphi_N$. Applying Lemma 1, we see that for all $n=0, 1,
\ldots, Q_N-1$
$$
\frac{1}{64(\Theta(N^N))^2}\leq
\int_{\mathbb T_N^3}e^{in\Phi_N} d\mu_{\mathbb T_N^3}.
$$
Therefore,
$$
\frac{1}{64(\Theta(N^N))^2}\leq
\int_{\mathbb T_N^3}\bigg(\frac{1}{Q_N}
\sum_{n=0}^{Q_N-1}e^{in\Phi_N}\bigg)
d\mu_{\mathbb T_N^3}.
$$
It remains to use identity (11). The lemma is proved.

\quad

   Our next aim is an upper bound for
$\delta_{\mathbb T_N^3, ~\mu_{\mathbb T_N^3}}(E_N)$. We shall
obtain it in Lemma 4. First we prove Lemma 3 which has a technical
character.

   Let $(\mathbb X_j, \mu_j), ~j=1, 2, \ldots m,$ be probability
spaces. let $(\mathbb X, \mu)$ be their product:
$$
\mathbb X=\mathbb X_1\times\mathbb X_2\times\ldots\times\mathbb X_m,
\qquad\mu=\mu_1\otimes\mu_2\otimes\ldots\otimes\mu_m.
$$
Consider a set $E\subseteq\mathbb X$. For each fixed $j=1, 2,
\ldots, m$ every point
$$
(x_1, x_2, \ldots, x_{j-1}, x_{j+1}, \ldots, x_m)
$$
in $\mathbb X_1\times\mathbb X_2\times\ldots \times\mathbb
X_{j-1}\times\mathbb X_{j+1}\times \ldots \times\mathbb X_m$
defines a $j$ -section of $E$, namely, the set
$$
E^{x_1, x_2,
\ldots, x_{j-1}, x_{j+1}, \ldots, x_m}=\{x_j\in\mathbb X_j :
(x_1, x_2, \ldots, x_{j-1}, x_j, x_{j+1}, \ldots, x_m)\in E\}.
$$

\quad

\textbf{Lemma 3.} \emph{Let $(\mathbb X_j, \mu_j), ~j=1, 2,
\ldots, m,$ be probability spaces and $(\mathbb X, \mu)$ their
product. Let $E$ be a $\mu$ -measurable set in $\mathbb X$. Let
$\delta\geq 0$. Suppose that for each $j=1, 2, \ldots, m$ every
$j$ -section $E^{x_1, x_2, \ldots, x_{j-1}, x_{j+1}, \ldots, x_m}$
of $E$ satisfies the condition that $\delta_{\mathbb X_j,
\mu_j}(E^{x_1, x_2, \ldots, x_{j-1}, x_{j+1}, \ldots, x_m})\leq
\delta$. Then $\delta_{\mathbb X, \mu}(E)\leq 3^{m-1}\delta$.}

\quad

\emph{Proof.}\footnote{Perhaps the assertion of Lemma 3 (in this
form or another) is known, but the author has not been able to
find an appropriate reference.} First we shall show that the
assertion of the lemma is true for $m=2$. Let $(\mathbb X_1,
\mu_1)$ and $(\mathbb X_2, \mu_2)$ be probability spaces. Put
$\mathbb X=\mathbb X_1\times\mathbb X_2$ and
$\mu=\mu_1\otimes\mu_2$. Consider an arbitrary $\mu$ -measurable
set $E\subseteq\mathbb X$ whose every $1$ -section
$$
E^{x_2}=\{x_1\in\mathbb X_1 : (x_1, x_2)\in E\}
$$
and every $2$ -section
$$
E^{x_1}=\{x_2\in\mathbb X_2 : (x_1, x_2)\in E\}
$$
satisfy
$$
\delta_{\mathbb X_1, \mu_1}(E^{x_2})\leq \delta
$$
and correspondingly
$$
\delta_{\mathbb X_2, \mu_2}(E^{x_1})\leq \delta.
$$
Let us verify that $\delta_{\mathbb X, \mu}(E)\leq 3\delta$.

   We put
$$
\mathbb X_1^<=\{x_1\in\mathbb X_1 : \mu_2(E^{x_1})\leq \delta\}, \qquad
\mathbb X_1^>=\{x_1\in\mathbb X_1 : \mu_2(E^{x_1})\geq 1-\delta\},
$$
and similarly
$$
\mathbb X_2^<=\{x_2\in\mathbb X_2 : \mu_1(E^{x_2})\leq \delta\}, \qquad
\mathbb X_2^>=\{x_2\in\mathbb X_2 : \mu_1(E^{x_2})\geq 1-\delta\}.
$$
We can assume that $0\leq\delta<1/2$ (for $\delta\geq 1/2$ the
assertion of the lemma is trivial). Then
$$
\mathbb X_j^<\cap\mathbb X_j^>=\varnothing, \quad
\mathbb X_j^<\cup\mathbb X_j^>=\mathbb X_j, \qquad j=1,2.
$$

   We put
$$
\alpha_1=\mu_1(\mathbb X_1^>), \qquad \alpha_2=\mu_2(\mathbb X_2^>).
$$

  Note that
$$
\iint_{\mathbb X_1^>\times\mathbb X_2^<}
1_E (x_1, x_2)d\mu_1(x_1)d\mu_2(x_2)\leq
\iint_{\mathbb X_1\times\mathbb X_2^<}
1_E (x_1, x_2)d\mu_1(x_1)d\mu_2(x_2)
$$
$$
=\int_{\mathbb X_2^<}\mu_1(E^{x_2})d\mu_2(x_2)
\leq \delta\mu_2(X_2^<)=\delta(1-\alpha_2).
\eqno(13)
$$
At the same time it is obvious that
$$
\iint_{\mathbb X_1^>\times\mathbb X_2^>}
1_E (x_1, x_2)d\mu_1(x_1)d\mu_2(x_2)\leq
\mu_1(\mathbb X_1^>)\mu_2(\mathbb X_2^>)=
\alpha_1\alpha_2.
\eqno(14)
$$
Adding (13) and (14), we obtain
$$
\iint_{\mathbb X_1^>\times\mathbb X_2}
1_E (x_1, x_2)d\mu_1(x_1)d\mu_2(x_2)\leq
\delta(1-\alpha_2)+\alpha_1\alpha_2.
\eqno(15)
$$
On the other hand
$$
\iint_{\mathbb X_1^>\times\mathbb X_2}
1_E (x_1, x_2)d\mu_1(x_1)d\mu_2(x_2)=
\int_{\mathbb X_1^>}\mu_2(E^{x_1})d\mu_1(x_1)
$$
$$
\geq (1-\delta)\mu_1(\mathbb X_1^>)=(1-\delta)\alpha_1.
\eqno(16)
$$
Together inequalities (15), (16) yield
$$
(1-\delta)\alpha_1\leq \delta(1-\alpha_2)+\alpha_1\alpha_2.
\eqno(17)
$$

   Similarly
$$
(1-\delta)\alpha_2\leq \delta(1-\alpha_1)+\alpha_2\alpha_1.
\eqno(18)
$$

  Adding inequalities (17) and (18), we see that
$$
\alpha_1+\alpha_2\leq 2\delta+2\alpha_1\alpha_2,
$$
that is
$$
\alpha_1(1-\alpha_2)+\alpha_2(1-\alpha_1)\leq 2\delta.
$$
Hence, putting $a=\min (\alpha_1, 1-\alpha_1)$, we have
$$
a=a(1-\alpha_2)+\alpha_2 a\leq
\alpha_1(1-\alpha_2)+\alpha_2(1-\alpha_1)\leq 2\delta.
$$
Thus, only two cases are possible: either $\alpha_1\leq 2\delta$
or $\alpha_1\geq 1-2\delta$.

   In the first case we obtain
\begin{multline*}
\mu(E)=\int_{\mathbb X_1}\mu_2(E^{x_1})d\mu_1(x_1)\\
=\int_{\mathbb X_1^>}\mu_2(E^{x_1})d\mu_1(x_1)+ \int_{\mathbb
X_1^<}\mu_2(E^{x_1})d\mu_1(x_1) \\ \leq \mu_1(\mathbb
X_1^>)+\delta\mu_1(\mathbb X_1^<)= \alpha_1+\delta(1-\alpha_1)\leq
\alpha_1+\delta\leq 3\delta.
\end{multline*}

  In the second case we obtain
\begin{multline*}
\mu(E)=\int_{\mathbb X_1}\mu_2(E^{x_1})d\mu_1(x_1)\geq
\int_{\mathbb X_1^>}\mu_2(E^{x_1})d\mu_1(x_1)\geq
(1-\delta)\mu_1(\mathbb X_1^>) \\ = (1-\delta)\alpha_1\geq
(1-\delta)(1-2\delta)= 1-3\delta+2\delta^2\geq 1-3\delta.
\end{multline*}
Thus, the assertion of the lemma is true for $m=2$.

   For an arbitrary $m$ the assertion of the lemma easily
follows by induction. Assume that the assertion of the lemma holds
for a certain $m\geq 2$. Let $(\mathbb X_j, \mu_j), ~j=1, 2,
\ldots, m+1,$ be probability spaces. It suffices to consider two
spaces, one of which is $\mathbb X_1\times\mathbb
X_2\times\ldots\times\mathbb X_m,$ with the measure
$\mu_1\otimes\mu_2\otimes\ldots\otimes\mu_m,$ and the other one is
$\mathbb X_{m+1}$ with the measure $\mu_{m+1}$. The lemma is
proved.

\quad

\textbf{Lemma 4.} \emph{If $N$ is a sufficiently large prime
number, then
$$
\delta_{\mathbb T_N^3, ~\mu_{\mathbb T_N^3}}(E_N)
\leq c\bigg(\frac{\log\log N}{\log N}\bigg)^{1/3}(\Theta(N^N))^2,
$$
where $c>0$ is independent of $N$.}

\quad

\emph{Proof.} Fix arbitrary $y\in\mathbb T_N$ and $z\in\mathbb
T_N$. Consider the corresponding $1$ -section of the set $E_N$,
namely, the set
$$
E_N^{y, z}=\{x\in\mathbb T_N : (x, y, z)\in E_N\}.
$$
Consider also the corresponding section $\Phi_N^{y, z}$ of the
function $\Phi_N$, i.e., the function of $x\in\mathbb T_N$,
defined as follows:
$$
\Phi_N^{y, z}(x)=\Phi_N(x, y, z).
$$
Note that in the product
$$
e^{in\Phi_N^{y, z}(x)}=
e^{in\varphi_N(x)}e^{in\varphi_N(z-x)}
e^{-in\varphi_N(y)}e^{-in\varphi_N(z-y)}
$$
only two factors are functions of $x$ both having $A(\mathbb
T_N)$ norm equal to $\|e^{in\varphi_N}\|_{A(\mathbb T_N)}$
while the two remaining factors are constants of modulus $1$.
So,
$$
\|e^{in\Phi_N^{y, z}}\|_{A(\mathbb T_N)}\leq
\|e^{in\varphi_N}\|_{A(\mathbb T_N)}^2.
$$
Hence, using Lemma 1, we see that
$$
\|e^{in\Phi_N^{y, z}}\|_{A(\mathbb T_N)}\leq (8\Theta(N^N))^2,
\qquad n=0, 1, 2, \ldots, Q_N-1.
$$
At the same time, due to identity (11), for the section $E_N^{y,
z}$ of the set $E_N$ we have
$$
1_{E_N^{y, z}}=\frac{1}{Q_N}\sum_{n=0}^{Q_N-1}e^{in\Phi_N^{y, z}}.
$$
Therefore,
$$
\|1_{E_N^{y, z}}\|_{A(\mathbb T_N)}\leq(8\Theta(N^N))^2.
$$

  Hence, using estimate (2), we see that for an
arbitrary $1$ -section $E_N^{y, z}$ of the set $E_N$ we have
$$
\delta_{\mathbb T_N, ~\mu_{\mathbb T_N}}(E_N^{y, z})
\leq c\bigg(\frac{\log\log N}{\log N}\bigg)^{1/3}(8\Theta(N^N))^2,
$$
where $c>0$ is independent of $y, z$ and $N$.

  In the same way we obtain similar upper estimate for an
arbitrary $2$ -section
$$
E_N^{x, z}=\{y\in\mathbb T_N : (x, y, z)\in E_N\}
$$
and an arbitrary $3$ -section
$$
E_N^{x, y}=\{z\in\mathbb T_N : (x, y, z)\in E_N\}.
$$

  Using Lemma 3 for $m=3$ and
$(\mathbb X_j, \mu_j)=(\mathbb T_N, \mu_{\mathbb T_N}), ~j=1, 2,
3,$ we see that
$$
\delta_{\mathbb T_N^3, ~\mu_{\mathbb T_N^3}}(E_N)
\leq 9 c\bigg(\frac{\log\log N}{\log N}\bigg)^{1/3}(8\Theta(N^N))^2.
$$
The lemma is proved.

\quad

  Consider the compliment $F_N=\mathbb T_N^3\setminus E_N$
of the set $E_N$.

\quad

\textbf{Lemma 5.} \emph{We have $\mu_{\mathbb
T_N^3}(F_N)\rightarrow 0$ as $N$ tends to $\infty$ running over
the set of prime numbers.}

\quad

\emph{Proof.} If $N$ is a sufficiently large prime number, then by
Lemma 4 we have either
$$
\mu_{\mathbb T_N^3}(E_N)\leq
c\bigg(\frac{\log\log N}{\log N}\bigg)^{1/3}(\Theta(N^N))^2
$$
or
$$
\mu_{\mathbb T_N^3}(F_N)\leq
c\bigg(\frac{\log\log N}{\log N}\bigg)^{1/3}(\Theta(N^N))^2.
\eqno(19)
$$

   In the first case, due to Lemma 2, we obtain
$$
\frac{1}{64(\Theta(N^N))^2}\leq
c\bigg(\frac{\log\log N}{\log N}\bigg)^{1/3}(\Theta(N^N))^2,
$$
that is
$$
\frac{1}{64c}\bigg(\frac{\log N}{\log \log N}\bigg)^{1/3}
\leq (\Theta(N^N))^4,
$$
which is impossible if $N$ is too large (see (3)). Thus, for all
sufficiently large prime $N$ we have estimate (19), whence, taking
(3) into account, we obtain
$$
\mu_{\mathbb T_N^3}(F_N)\leq c\bigg(\frac{\log\log N}{\log N}\bigg)^{1/3}
o\bigg(\bigg(\frac{\log N}{\log\log N}\bigg)^{1/6}\bigg)=o(1).
$$
The lemma is proved.

\quad

  Now we shall conclude the proof of the theorem. Note that
(since $e^{i\Phi_N}=1$ on $E_N$) estimate (7) yields
\begin{multline*}
\|e^{i\Phi}-1\|_{L^\infty(E_N)}=
\|e^{i\Phi}-e^{i\Phi_N}\|_{L^\infty(E_N)} \leq
\|\Phi-\Phi_N\|_{L^\infty(E_N)} \\ \leq
\|\Phi-\Phi_N\|_{L^\infty(\mathbb T_N^3)}\leq
 4 \|\varphi-\varphi_N\|_{L^\infty(\mathbb T_N)}
\leq\frac{8\pi}{NQ_N}\leq\frac{8\pi}{N}=o(1), \qquad
N\rightarrow\infty.
\end{multline*}

  Hence, using Lemma 5, we see that as $N$
tends to infinity running over prime numbers, we have
\begin{multline*}
\int_{\mathbb T_N^3} |e^{i\Phi}-1| d\mu_{\mathbb T_N^3}=
\int_{E_N} |e^{i\Phi}-1| d\mu_{\mathbb T_N^3}+ \int_{F_N}
|e^{i\Phi}-1| d\mu_{\mathbb T_N^3}\\ \leq
 \|e^{i\Phi}-1\|_{L^\infty(E_N)}+ 2\mu_{\mathbb
T_N^3} (F_N) \rightarrow 0.
\end{multline*}

  At the same time (since the function $\Phi$ is continuous
on $\mathbb T^3$), we have
\begin{multline*}
\int_{\mathbb T_N^3} |e^{i\Phi}-1| d\mu_{\mathbb T_N^3} \\ =
\frac{1}{N^3}\sum_{0\leq k, l, m\leq N-1} |e^{i\Phi(\frac{2\pi
k}{N}, \frac{2\pi l}{N}, \frac{2\pi m}{N})}-1|\rightarrow
\frac{1}{(2\pi)^3}\int_{\mathbb T^3}|e^{i\Phi(x, y, z)}-1| dxdydz
\end{multline*}
as $N\rightarrow\infty$.
  Thus,
$$
\int_{\mathbb T^3}|e^{i\Phi(x, y, z)}-1| dxdydz=0.
$$

  Hence (since $\Phi$ is continuous) we obtain
$$
e^{i\Phi(x, y, z)}=1
$$
for all $(x, y, z)\in\mathbb T^3$ and we see that
$$
\Phi(x, y, z)=2\pi k, \qquad (x, y, z) \in \mathbb T^3,
$$
where $k$ is an integer independent of $(x, y, z)$. Taking $(x, y,
z)=(0, 0, 0)$ we obtain $k=0$, whence $\Phi(x, y, z)=0$ on
$\mathbb T^3$. That is
$$
\varphi(x)+\varphi(z-x)-\varphi(y)-\varphi(z-y)=0
$$
for all $x, y, z \in \mathbb R$. The assertion of the theorem
follows.

\quad

\textbf{Remarks.}

   1. The theorem obtained in this work has the following
(equivalent) operator version. Let $U$ be a bounded translation
invariant operator from $l^1(\mathbb Z)$ to itself such that
$$
\|U^n\|_{l^1\rightarrow l^1}=
o\bigg(\bigg(\frac{\log\log |n|}{\log\log\log |n|}\bigg)^{1/12}\bigg),
\qquad |n|\rightarrow\infty, \quad n\in \mathbb Z,
\eqno(20)
$$
then $U=\xi S$, where $\xi$ is a complex number, $|\xi|=1$, and
$S$ is a translation. Indeed, it is easy to verify (and is well
known, see, e.g., [13]) that each bounded translation invariant
operator from $l^1(\mathbb Z)$ to itself is an operator of
convolution with a certain sequence that belongs to $l^1(\mathbb
Z)$ and the norm of the operator is equal to the $l^1$ norm of the
sequence. In particular, $U$ is an operator of convolution with a
certain sequence $\{u_k, ~k\in \mathbb Z\}\in l^1$. Define a
function $u$ on $\mathbb T$ by
$$
u(t)=\sum_{k\in\mathbb Z} u_k e^{ikt}.
$$
We have $u\in A(\mathbb T)$ and
$$
\|u^n\|_{A(\mathbb T)}=\|U^n\|_{l^1\rightarrow l^1}.
$$
The function $u$ is continuous and since
$\|u^n\|_{L^\infty(\mathbb T)}\leq\|u^n\|_{A(\mathbb T)}$, it is
clear that $|u(t)|=1$ for all $t\in\mathbb T$. (Otherwise the
growth of the norms $\|U^n\|_{l^1\rightarrow l^1}$ would be
exponential either as $n\rightarrow +\infty$ or as $n\rightarrow
-\infty$.) Thus $u$ is a continuous function that maps $\mathbb R$
into the circle $\mathcal{C}=\{z\in \mathbb C : |z|=1\}$ on the
complex plane $\mathbb C$. Every such function has the form
$u(t)=e^{i\varphi(t)}$ where $\varphi : \mathbb
R\rightarrow\mathbb R$ is a continuous function. \footnote{The
function $\varphi$ is the lifting of $u$ with respect to the
covering $e: \mathbb R\rightarrow \mathcal{C}$ where
$e(x)=e^{ix}$. The existence of $\varphi$ is guaranteed by the
monodromy theorem, see, e.g., [14, Lecture 4, corollary from
Theorem 1].} Since the function $u$ is $2\pi$ -periodic, we have
$\varphi(t+2\pi)=\varphi(t)(\mathrm{mod} ~2\pi)$, i.e., $\varphi$
is a continuous mapping of the circle $\mathbb T$ into itself. We
have $\|e^{in\varphi}\|_{A(\mathbb T)}=\|u^n\|_{A(\mathbb
T)}=\|U^n\|_{l^1\rightarrow l^1}$. Using our theorem, we obtain
from (20) that $\varphi(t)=\nu t+\varphi(0)$ for some
$\nu\in\mathbb Z$. So $u(t)=\xi e^{i\nu t}$, where
$\xi=e^{i\varphi(0)}$. Therefore, the operator $U$ is the operator
of convolution with the sequence $\{u_k, ~k\in\mathbb Z\}$, where
$u_k=0$ for all $k\neq \nu$ and $u_\nu=\xi$. Thus, $U=\xi S$,
where $\xi\in\mathbb C, ~|\xi|=1,$ and $S$ is a translation ($S :
\{x_k\}\rightarrow \{x_{k-\nu}\}$).

  2. The theorem of the present work easily transfers to
the multidimensional case. Let $A(\mathbb T^d)$ be the space of
all continuous functions on the $d$ -dimensional torus $\mathbb
T^d$ with absolutely convergent Fourier series. We put
$\|f\|_{A(\mathbb T^d)}=\|\widehat{f}\|_{l^1(\mathbb Z^d)}$, where
$\widehat{f}$ is the Fourier transform on $\mathbb T^d$. If
$\varphi$ is a continuous mapping $\mathbb T^d\rightarrow\mathbb
T$ such that
$$
\|e^{in\varphi}\|_{A(\mathbb T^d)}=
o\bigg(\bigg(\frac{\log\log |n|}{\log\log\log |n|}\bigg)^{1/12}\bigg),
\qquad n\in\mathbb Z, \quad |n|\rightarrow \infty,
$$
then $\varphi$ is linear, i.e., $\varphi(t)=(\nu,t)+\varphi(0)$,
where $\nu\in\mathbb Z^d$. (Here $(\nu, t)$ stands for the inner
product of vectors $\nu\in\mathbb Z^d$ and $t\in\mathbb T^d$.) The
multidimensional case easily reduces to the one-dimensional case
by induction over the dimension. It suffices only to note the
following. Let $f\in A(\mathbb T^{d+1})$. For a fixed $x\in\mathbb
T^d$ consider the function $f_x$ on $\mathbb T$ defined by
$f_x(y)=f(x, y)$. Then $\|f_x\|_{A(\mathbb T)}\leq\|f\|_{A(\mathbb
T^{d+1})}$.

   The corresponding version for operators $U$ acting in
$l^1(\mathbb Z^d)$ also holds.

  3. As is noted by B. Green and S. V. Konyagin
[10, \S~1], it is plausible that the estimate in their theorem can
be improved. Perhaps $(\log\log N/\log N)^{1/3}$ on the right side
in (2) can be replaced with $1/\log N$ (further improvement is
impossible). This improvement would allow to replace the right
side in the condition of our theorem with $o((\log\log
|n|)^{1/4})$. Apparently it is easier to obtain an improvement of
estimate (2) with the replacement of the exponent $1/3$ with
$1/2$; in this connection see the work of T. Sanders [15]. This
would allow to replace the exponent $1/12$ in our theorem with
$1/8$. As to the possibility to eliminate $\log\log N$ in the
Green--Konyagin theorem, see [10, \S~6]. This would allow to
eliminate $\log\log\log |n|$ in our theorem.

  4. Kahane's conjecture remains unproved even under the additional
assumption of $C^1$ -smoothness of $\varphi$. In this connection
note the following result of the author [6] (see also [7]): if
$\gamma(n)\geq 0$ is an arbitrary sequence,
$\gamma(n)\rightarrow+\infty$, then there exists a $C^1$ -smooth
non-linear mapping $\varphi : \mathbb T\rightarrow\mathbb T$ such
that $\|e^{in\varphi}\|_{A(\mathbb T)}=O(\gamma(|n|)\log |n|)$.

  5. It is not clear if there exists a continuous non-linear
mapping $\varphi : \mathbb T\rightarrow\mathbb T$ satisfying
$$
\|e^{in_k\varphi}\|_{A(\mathbb T)}=O(1),
\eqno(21)
$$
for some unbounded sequence of integers $\{n_k\}$. If $\varphi$ is
absolutely continuous (i.e., if the function $\varphi: \mathbb
R\rightarrow\mathbb R$ with condition (1) is absolutely continuous
on each interval of length  $2\pi$) then this is impossible. We
shortly explained this in [6]. Here we shall give a detailed and
simple proof.

   Note first that if $g$ is a real measurable function on
$\mathbb T$ such that $\widehat{e^{i\lambda g}}\in l^1(\mathbb Z)$
for all $\lambda\in\mathbb R$ and
$$
\sup_{\lambda\in\mathbb R}
\|\widehat{e^{i\lambda g}}\|_{l^1(\mathbb Z)}<\infty,
$$
then $g$ is constant almost everywhere. One can see this as
follows. We have $\widehat{e^{ig}}\in l^1(\mathbb Z)$, therefore,
$e^{ig}$ coincides almost everywhere with a certain continuous
function $\xi$. Then $|\xi|=1$ almost everywhere and therefore
everywhere. It follows that $\xi=e^{i\psi}$, where $\psi : \mathbb
T\rightarrow\mathbb T$ is a continuous mapping (this has already
been explained in the first remark). Then for each $n\in\mathbb Z$
we have $e^{in\psi}=\xi^n=e^{ing}$ almost everywhere, therefore,
$\|e^{in\psi}\|_{A(\mathbb T)}=\|\widehat{e^{ing}}\|_{l^1(\mathbb
Z)}=O(1), ~n\in\mathbb Z$. By the Beurling--Helson theorem the
function $\psi$ is linear with integer slope. Thus, $g$ almost
everywhere coincides $\mathrm{mod}~2\pi$ with a linear function
that has an integer slope. Let us take a (real) irrational number
$\alpha$. Repeating the arguments for the function
$g_\alpha=\alpha g$ instead of $g$, we see that $\alpha g$ is also
almost everywhere coincides $\mathrm{mod}~2\pi$ with a linear
function that has an integer slope. This is possible only in the
case when $g$ is constant almost everywhere.

  Let $\varphi : \mathbb T\rightarrow\mathbb T$ be an absolutely
continuous mapping. Suppose that we have (21). Subtracting an
appropriate linear function we can assume that $\varphi$ is a real
$2\pi$ -periodic function on $\mathbb R$. We can also assume that
$n_k\rightarrow\infty$. For arbitrary $\lambda\in\mathbb R,
~\lambda\neq 0,$ and $k=1, 2, \ldots$ define a function
$g_{\lambda, k}$ on $\mathbb R$ by
$$
g_{\lambda, k}(t)=
\frac{\varphi(t+\lambda/n_k)-\varphi(t)}{\lambda/n_k}.
$$
As $k\rightarrow\infty$, we have $g_{\lambda, k}\rightarrow
\varphi'$ almost everywhere. Note that
$$
e^{i\lambda g_{\lambda, k}(t)}=
e^{in_k\varphi(t+\lambda/n_k)}e^{-in_k\varphi(t)},
$$
whence, assuming that $\|e^{in_k\varphi}\|_{A(\mathbb T)}\leq c,
~k=1, 2, \ldots,$ we obtain $\|e^{i\lambda g_{\lambda,
k}}\|_{A(\mathbb T)}\leq c^2$. Tending $k$ to $\infty$, we see
that $\|\widehat{e^{i\lambda\varphi'}}\|_{l^1(\mathbb Z)}\leq c^2$
for all $\lambda\in\mathbb R$ (the case when $\lambda=0$ is
obvious). Thus, the derivative $\varphi'$ of the function
$\varphi$ is constant almost everywhere. From the condition of
absolute continuity it follows that $\varphi$ is linear.

\quad

  Initially, the theorem obtained in the present work
was somewhat weaker. I am grateful to S. V. Konyagin who turned my
attention to Theorem 1.3 of the work [10]. This allowed to improve
the result.

  I am also grateful to Yu. N. Kuznetsova for the help with the
proof of the lemma on the sections (Lemma 3).

\quad

\begin{center}
\textbf{References}
\end{center}
\flushleft
\begin{enumerate}

\item A. Beurling, H. Helson, ``Fourier-Stieltjes transforms
    with bounded powers'', \emph{Math. Scand.,} \textbf{1}
    (1953), 120-126.

\item J.-P. Kahane, \emph{S\'erie de Fourier absolument
    convergentes}, Springer-Verlag, Berlin--Heidelberg--New
    York, 1970.

\item J.-P. Kahane, ``Quatre le\c cons sur les
    hom\'eomorphismes du circle et les s\'eries de Fourier'',
    in: \emph{Topics in Modern Harmonic Analysis,} Vol. II,
    Ist. Naz. Alta Mat. Francesco Severi, Roma, 1983, 955-990.

\item  Z. L. Leibenson, ``On the ring of functions with
    absolutely convergent Fourier series'', \emph{Uspehi
    Matem. Nauk}, \textbf{9}:3(61) (1954), 157-162 (in
    Russian).

\item J.-P. Kahane, ``Sur certaines classes de s\'eries de
    Fourier absolument convergentes'', \emph{J. de
    Math\'ematiques Pures et Appliqu\'ees}, \textbf{35}:3
    (1956), 249-259.

\item  V. V. Lebedev, ``Diffeomorphisms of the circle and the
    Beurling--Helson theorem'', \emph{Functional analysis and
    its applications}, \textbf{36}:1(2002), 25-29.

\item V. V. Lebedev, ``Quantitative estimates in
    Beurling--Helson type theorems'', \emph{Sbornik:
    Mathematics}, \textbf{201}:12 (2010), 1811-1836.

\item V. V. Lebedev, ``Estimates in Beurling--Helson type
    theorems: Multidimensional case'', \emph{Mathematical
    Notes}, \textbf{90}:3 (2011), 373-384.

\item J.-P. Kahane, ``Transform\'ees de Fourier des fonctions
    sommables'', \emph{Proceedings of the Int. Congr. Math.,
    15-22 Aug., 1962, Stockholm, Sweden}, Inst.
    Mittag-Leffler, Djursholm, Sweden, 1963, pp. 114-131.

\item B. Green, S. Konyagin, ``On the Littlewood problem
    modulo a prime'', \emph{Canad. J. Math.}, \textbf{61}:1
    (2009), 141-164.

\item E. M. Stein and R. Shakarchi, \emph{Fourier analysis: An
    introduction (Princeton Lectures in Analysis v. I)},
    Princeton Univ. Press, Princeton and Oxford, 2003.

\item W. M. Schmidt, \emph{Diophantine Approximation}, Lect.
    Notes in Math. 785, Springer-Verlag, Berlin-Heidelberg-New
    York, 1980.

\item R. Larsen, \emph{An introduction to the theory of
    multipliers}, Springer-Verlag, Berlin-Heidelberg-New York,
    1971.

\item M. M. Postnikov, \emph{Le\c cons de g\'eom\'etrie.
    Semester IV. G\'eom\'etrie diff\'erentielle}, Mir, Moscow,
    2001.

\item T. Sanders, ``The Littlewood--Gowers problem'',
    \emph{Journal d' Analyse Mathematique}, \textbf{101}:1
    (2007), 123-162.

\end{enumerate}

\quad

\qquad \textsc{V. V. Lebedev}\\
\qquad Dept. of Mathematical Analysis\\
\qquad Moscow State Institute of Electronics\\
\qquad and Mathematics (Technical University)\\
\qquad E-mail address: \emph {lebedevhome@gmail.com}

\end{document}